\documentclass[12pt]{amsart}

\usepackage{graphicx, verbatim,amsmath,amssymb, float}
\usepackage{hyperref}
\hypersetup{
    colorlinks=true,
    linkcolor=blue,
    filecolor=magenta,      
    urlcolor=cyan,
}

\def\vs#1 {\vskip#1truein}
\def\hs#1 {\hskip#1truein}

\newcommand{\be}{\begin{equation}}
\newcommand{\ee}{\end{equation}}

\newcommand{\Z}{{\mathbb Z}}

\newtheorem*{thm*}{Theorem}

\everymath{\displaystyle}

\theoremstyle{definition}

\title{Conway and Aperiodic Tilings}

\date{\today}

\begin{document}


\maketitle
\centerline{Charles Radin}
\vs.1 
\centerline{Department of Mathematics}

\centerline{The University of Texas at Austin}

\centerline{Austin, TX 78712 }
\section{The origins of aperiodicity}\label{aperiodicity}
\begin{figure}[H]
\center{\includegraphics[width=7in]{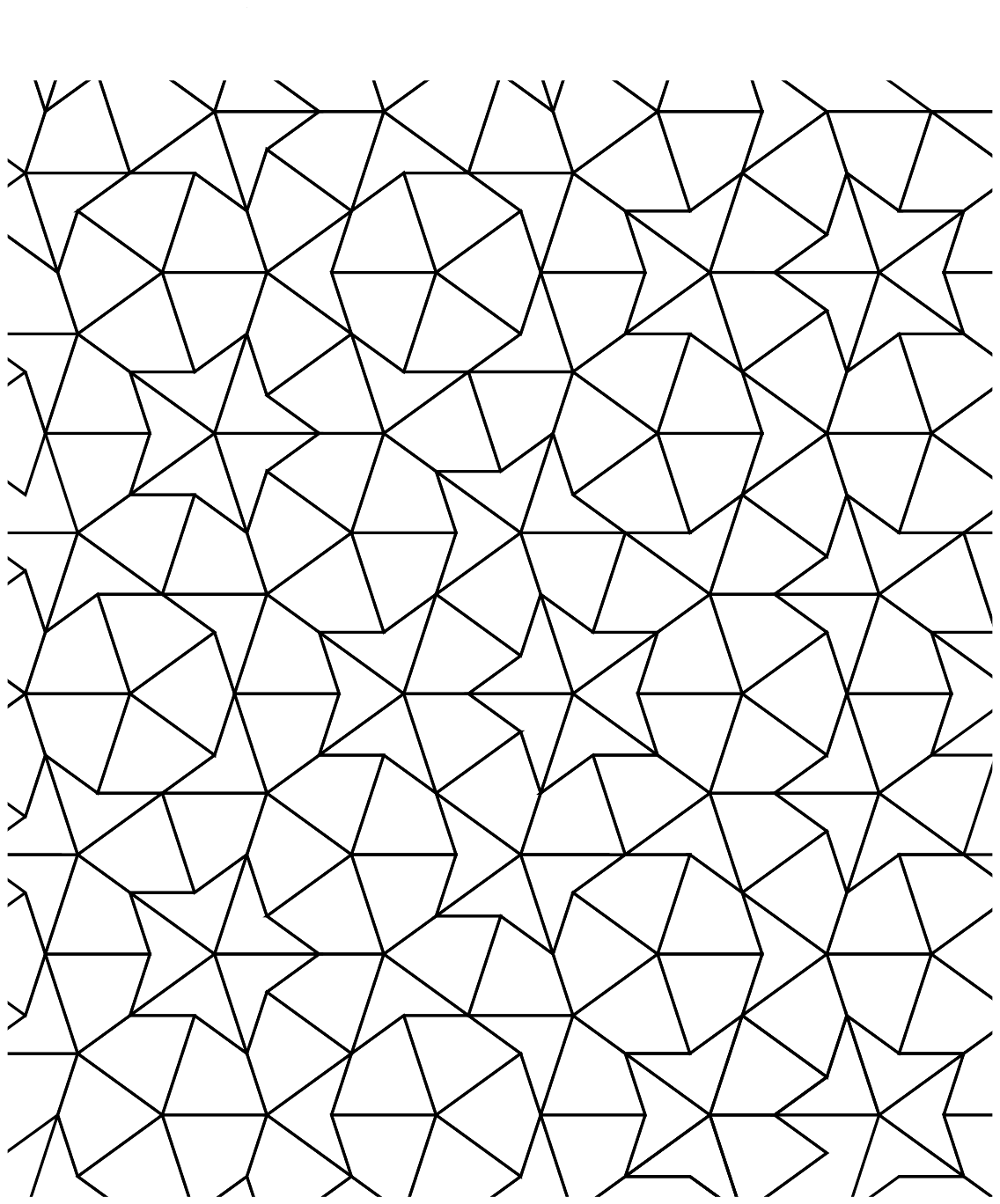}}
\caption{Kite and Dart tiling.}
\label{KandDtiling}
\end{figure}

John Conway became interested in the kite and dart tilings (Fig.\ ~\ref{KandDtiling})
soon after their discovery by Roger Penrose around 1976, and was well
known for his analysis of them and the memorable terminology he
introduced. (John has enhanced many subjects with colorful
terminology, probably any subject on which he spent much time.) Neither
Penrose nor Conway actually published much about the tilings; word
spread initially through talks given by John, at least until Martin
Gardner wrote a wonderful Mathematical Games article in Scientific
American (1977) \cite{Gar}, tracing their history back to mathematically simpler
tilings by Robert Berger, Raphael Robinson and others. I stumbled on
the Gardner article in 1983 and was interested in the earlier tilings,
on realizing they could be used to make unusual  models
of solids in statistical physics.
I worked on this connection for a number of years, and when Conway was
visiting the knot theory group in Austin in 1991 I asked him about a
technical point on which I had been thinking. Within two hours he
answered my question, and for the next 10 years or so we had frequent
correspondence, which among other things resulted in three papers. It
was great fun working with John and I learned a lot, in part due to
his unusual approach to mathematical research. I will describe our
interaction and the results that came out of it in that decade, but to
place the work in context I will start with the prehistory as I first
learned it from Gardner's article.

In 1961 the logician Hao Wang considered the following
structure \cite{Wan}. Start with a finite number of unit squares in the plane,
and deform their edges with small bumps and dents.  Wang was interested in
tilings of the plane by translations of these `tiles', a
tiling being a covering without overlap. (Think of the squarish tiles
in a tiling fitting together as in an infinite jigsaw puzzle, bumps
fitting into neighboring dents and with corners touching and producing a copy of
$\Z^2$.)

Wang asked if there existed such a finite collection of tiles (i.e.\
of edge deformations of unit squares) with these two properties:
\begin{itemize}
\item[a)] one could tile the plane by {\em translates} of the tiles;
\item[b)] {\em
  no} such tiling is periodic, that is, invariant under a pair of
linearly independent translations.
\end{itemize}
Intuitively he was asking if there was a finite collection of deformed
squares such that translation copies could tile the plane but {\em
  only} in complicated ways, namely nonperiodically.  (Other notions
of `complicated', such as `algorithmically complex', were considered by
others later, but this is the original, brilliant idea.)

\begin{figure}[H]
\center{\includegraphics[width=6in]{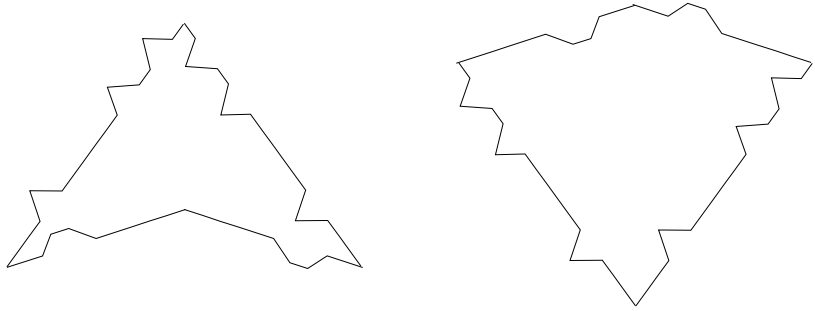}}
\caption{Kite and Dart tiles.}
\label{KandDtiles}
\end{figure}

Wang was interested in this question because the solution would answer
an open decidability problem in predicate calculus, an indication
of its significance. His student Berger finally solved the tiling
problem in the affirmative in his thesis (1966) \cite{Ber}, though
Wang had already solved the logic problem another way. Berger found a
number of such `aperiodic tiling sets', well motivated but all rather
complicated. Over the years nicer solutions were found, in particular
by Robinson, when the example by Penrose appeared in Gardner's column
in 1977: two (nonsquare) shapes, a kite and a dart, and after
deforming their edges (not shown in the tiling) they appear in Fig.\
~\ref{KandDtiles}.

This example doesn't quite fit Wang's structure though the
difficulties are easily resolved. The fact that the tiles are not
square isn't so important since they are still compatible with
nonorthogonal axes. And since the tiles only appear in a tiling in 10
different orientations one could simply call rotated tiles different
tile types, and not need rotations. In any case the tilings were
certainly interesting to many readers of Gardner's column.

Wang's invention of `aperiodicity' evolved in several directions,
including physics modeling, which was my original interest in
them. This article celebrates John Conway and his mathematics and we
will concentrate on the {\em geometric} developments coming from
the kite and dart tilings, rather than the broader subject.  

\section{The pinwheel and quaquaversal tilings}\label{pinwheel}
\begin{figure}[H]
\center{\includegraphics[width=5in]{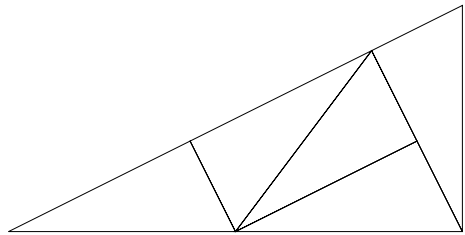}}
\caption{Pinwheel substitution.}
\label{pinwheel_substitution}
\end{figure}

When John visited Austin in 1991 I pointed out that in all known
aperiodic tiling sets, for instance the kite and dart, it was a
curious fact that the tiles in each such tiling only appear in a
finite number of relative orientations. By analogy to particle
configurations this seemed artificial and I thought there ought to be
examples without that special feature. On the face of it this was a
hard problem -- any new aperiodic tiling set is hard to find, even
without this new property -- but there was a natural path which I
thought might yield such an example if only one could solve a certain
much simpler problem first. I asked John if there was some polygon in
the plane such that it could be decomposed into smaller copies of
itself, all the same size, so that their relative orientations
included an angle irrational with respect to $\pi$. He quickly came up
with the $1,2,\sqrt{5}$ right triangle with the decomposition in Fig.\
~\ref{pinwheel_substitution}.  This did not solve my aperiodic tiling
problem, but I hoped I could now solve that using a beautiful recent
result of Shahar Mozes concerning tilings with squarish tiles.  
In his Master's thesis (1989) \cite{Moz} Mozes had studied the technique of
Berger and Robinson and broadly generalized it.  I will give one
ingredient in Mozes' theorem through an example.

A simple method for building a complicated sequence of symbols is by
substitution. For instance the Thue-Morse sequence 
$$abbabaabbaababba\ldots $$
can be generated by
repeating the substitution: $a\to ab,\, b\to ba$, obtaining
consecutively:
$$a \to ab \to abba \to abbabaab \to  \ldots$$

\begin{figure}[H]
\center{\includegraphics[width=5in]{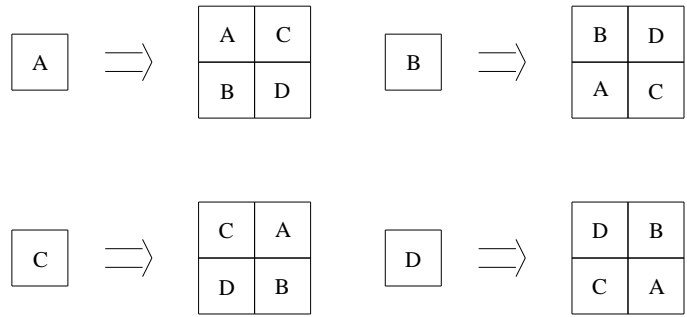}}
\caption{Morse substitution.}
\label{Morse_substitution}
\end{figure}

\begin{figure}[H]
\center{\includegraphics[width=5in]{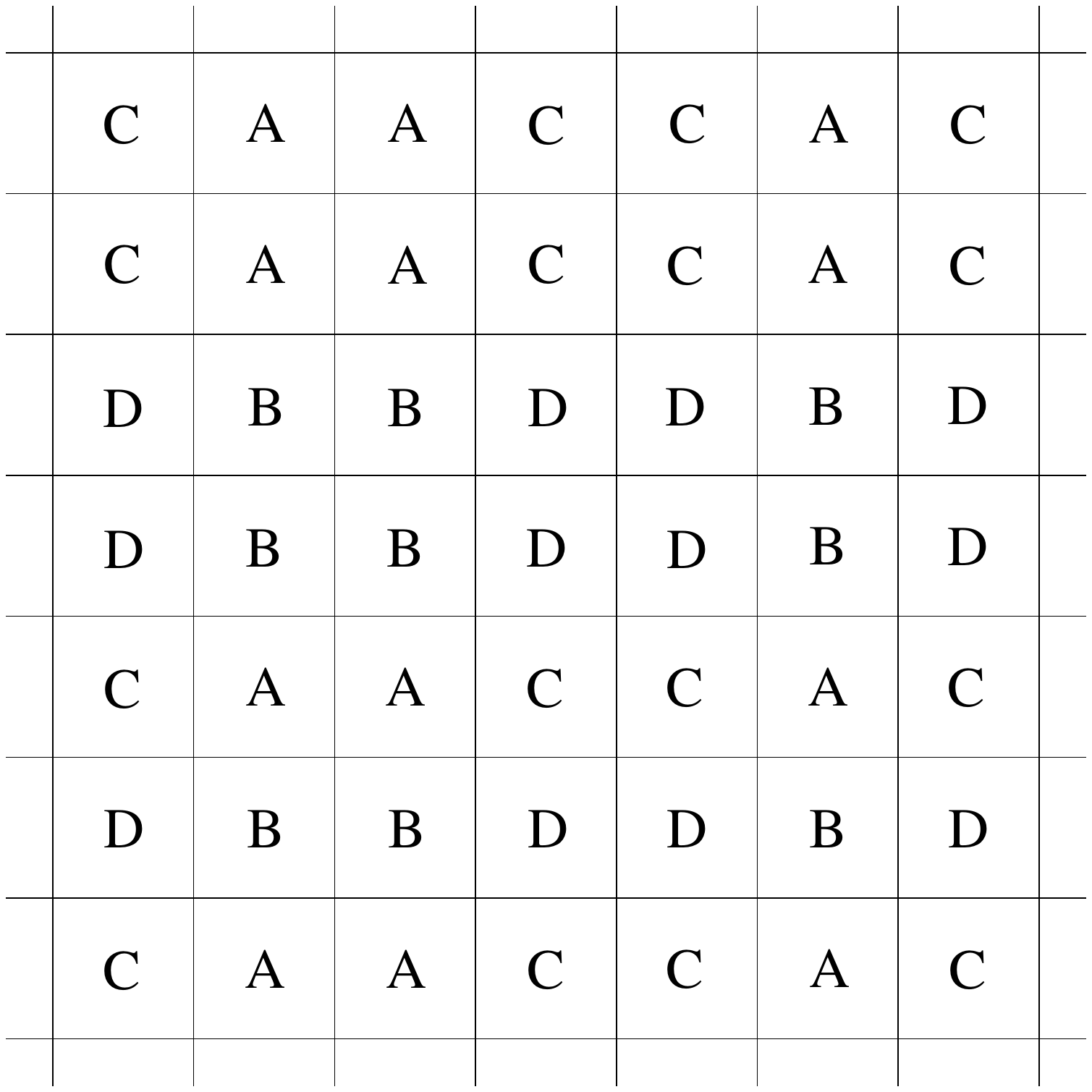}}
\caption{Morse tiling.}
\label{Morse_tiling}
\end{figure}

\noindent Applying this idea to arrays in the plane, the substitution in Fig.\
~\ref{Morse_substitution} gives rise to the infinite array in Fig.\
~\ref{Morse_tiling}. (To get many noncongruent tilings of the whole
plane one can `expand' about different points each time.)  

Mozes used substitutions to generalize the technique by Robinson to
add bumps and dents to the edges of the unit squares so that the
tilings by the jigsaw pieces are precise analogues of those in the
infinite substitutions.  (In his language he was studying subshifts
with $\Z^2$ action, and proved how, given a substitution subshift,
to construct a subshift of finite type which was conjugate to it as a
dynamical system.)  It should be noted that this process usually
produces many differently deformed copies of each of the original
square shapes.

\begin{figure}[H]
\center{\includegraphics[width=5in]{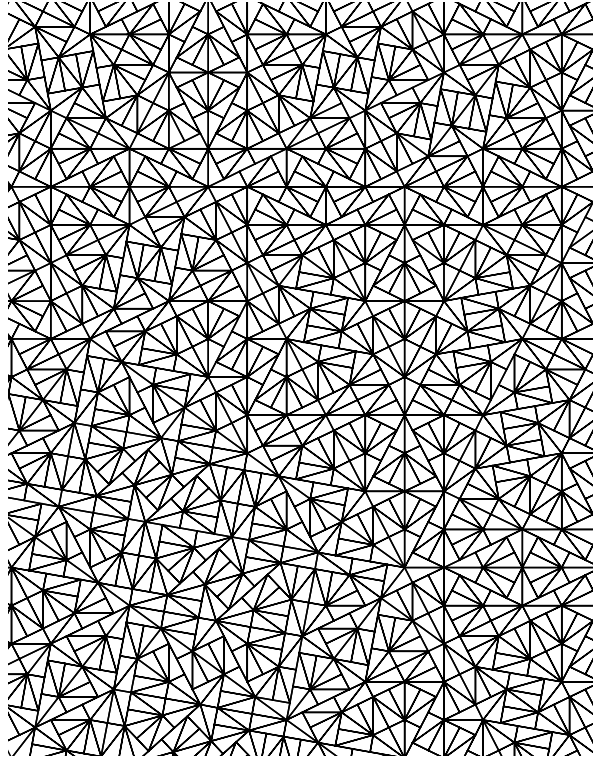}}
\caption{Pinwheel tiling.}
\label{pinwheel_tiling}
\end{figure}

Iterating  the pinwheel substitution of Fig.\ ~\ref{pinwheel_substitution} one gets 
substitution pinwheel tilings such as in Fig.\ ~\ref{pinwheel_tiling}, which have
the desired range of orientations.

Before going further I should mention a complication, in that if I 
succeeded in finding deformed triangles which only tile the plane
as in Fig.\ ~\ref{pinwheel_tiling}, it would require infinitely many 
different tiles if one only used translations as in a). What I was hoping
to do was find a {\em finite} number of tiles such that:
\begin{itemize}
\item[a')] {\em congruent } copies can tile the plane;
\item[b)] no such
tiling is periodic.
\end{itemize}

Note that geometry is slipping in by the replacement of the
translation symmetries of $\Z^2$ by congruences of the plane in a'), a
change which is not in the spirit  of Wang's
invention, which was purely algorithmic, but which opened the door
to interesting geometry, as we will see.

The square grid  is heavily used by Mozes (and Berger and Robinson), so it was not
at all clear that the technique could be adapted when the triangular
tiles are also rotated as in the substitution of Fig.\ ~\ref{pinwheel_substitution}. But I eventually
adapted it (1992), though the proof is much more complicated than that
of Mozes.  This gave a large but finite set of differently deformed
triangular tiles which only tile as in Fig.\ ~\ref{pinwheel_tiling}.
John was not interested in this. I
worked with him for years and it was clear from the beginning that he
was very particular in what he would be willing to think about.  What
attracted him to the kite and dart tilings was the challenge, which
motivated Penrose also, to find an aperiodic tiling set with as few
tiles as possible.  He loved puzzles, in particular about geometry in
the sense of Euclid.  He loved discrete math in various senses.  What
he did {\em not} love was analysis, and there was no way I could get
him interested in the ergodic theory of Mozes, work which I found beautiful.

Anyway, the pinwheel example brought into play the rotations in the
plane in making aperiodic tiling sets, and is mathematically a shift
in the subject towards geometry, in the sense of Klein's
Erlangen program, geometry viewed through the isometry groups of spaces.

Since the rotation group in three dimensions is much richer than in
two dimensions, it was natural to see what would happen if one started
with a substitution there. It was not clear that I could interest John
in such a project. Indeed, once my very complicated proof of the
pinwheel bumps and dents (more formally called `matching rules') was
published, John sat me down in his office and informed me that I had
wasted a lot of effort and that he would show me how to obtain
matching rules for the pinwheel substitution in a very simple way. The
idea was to look inside a pinwheel tiling at the `neighborhoods' of
tiles, made not just by nearest neighbors, but next nearest neighbors
etc., to some appropriate fixed finite depth, and use these to define
tiles with matching rules. This technique in fact was known to work
for some older examples coming from a substitution and John was sure
it was true generally. So we spent some time on the floor looking at
larger and larger neighborhoods in a pinwheel tiling drawn on a large
piece of paper. Finally John let me show him a proof by Ludwig Danzer
that this method cannot work for the pinwheel. John was {\em very}
surprised. And he was now sufficiently appreciative of what I did on
the pinwheel that our discussions became {\em almost} 
two-sided. 

\begin{figure}[H]
\center{\includegraphics[width=5in]{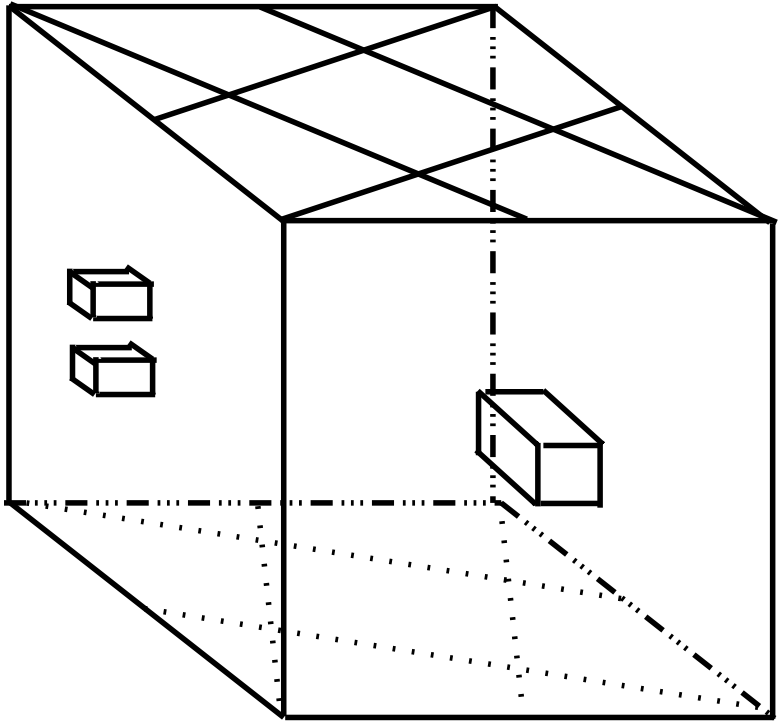}}
\caption{Schmitt tile.}
\label{Schmitt_tile}
\end{figure}

Also about this time I told him of an old but unpublished example
(1988) by Peter Schmitt, a
single deformed unit cube in space, which was aperiodic. It was based
on the following simple idea. One puts bumps and dents on four sides
of the cube so it can only tile space if it forms infinite slabs of unit thickness:
think of cubes lined up on a table, checkerboard fashion.  Then on the
tops cut two or three straight grooves (i.e.\ dents) at an appropriate angle with respect
to the edges, and at the bottom corresponding raised bumps (the reverse of the
dents), but at an angle with respect to the dents which is irrational
with respect to $\pi$. (See Fig.\ ~\ref{Schmitt_tile}.) If you think for a moment it is clear that this
will only allow tilings which consist of infinite parallel slabs
rotated irrationally with respect to one another, and therefore
nonperiodic. Wonderfully simple. Anyway, John was {\em very}
interested. He had gotten interested in the kite and dart tilings as
the smallest known aperiodic set -- two tiles -- and had spent a lot
of time trying to find a single planar tile which was aperiodic,
without success. And here was a simple way to find such a tile in
three dimensions! Immediately he wanted to improve it. We spent a few
hours (again on the floor, over a large sheet of paper) smoothing out
the bumps and dents, and came up with what was arguably a nicer
example. This illustrates well our different interests in
aperiodicity. For John it was the challenge of finding simpler and
nicer examples of aperiodic tiling sets, and although I tried I never
understood the point.

\begin{figure}[H]
\center{\includegraphics[width=5in]{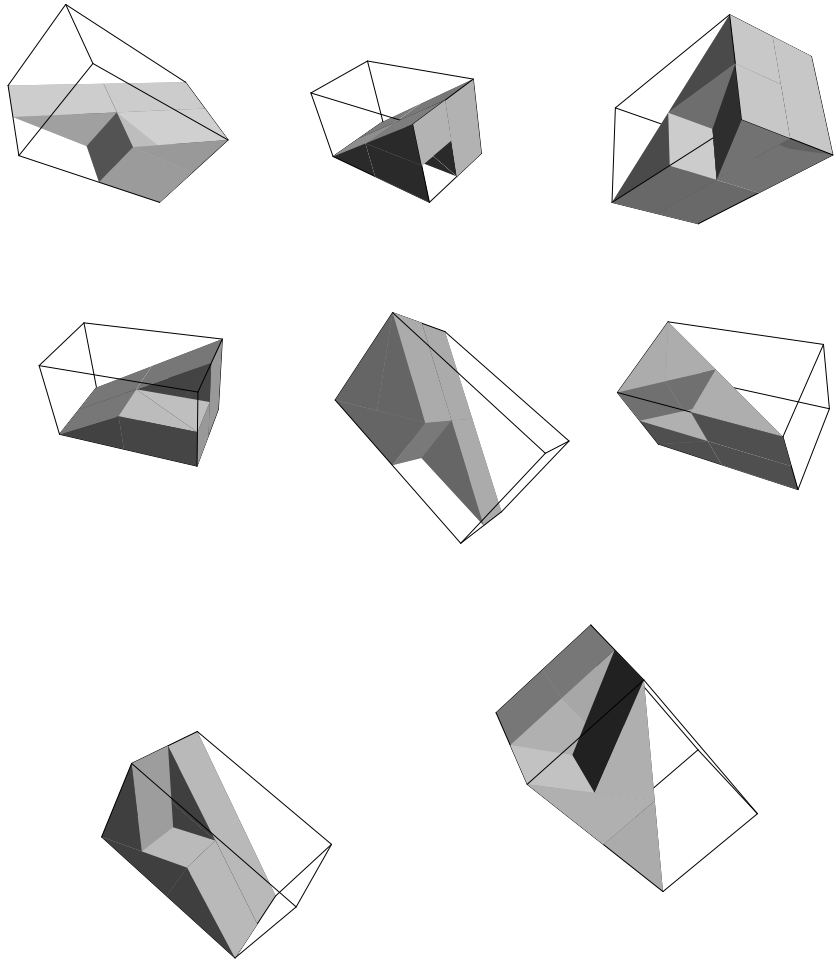}}
\caption{Quaquaversal substitution.}
\label{Quaquaversal_substitution}
\end{figure}

\begin{figure}[H]
\center{\includegraphics[width=5in]{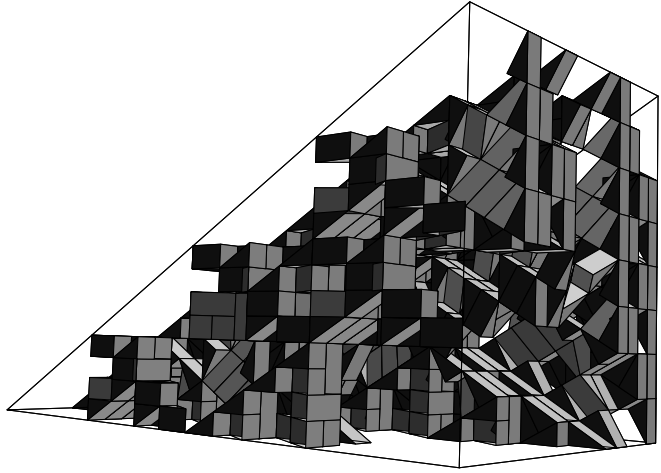}}
\caption{Quaquaversal tiling.}
\label{Quaquaversal_tiling}
\end{figure}

In 1995 John was willing to consider my
proposal to investigate tilings in three dimensions exhibiting special
properties.  Rather quickly we
(mostly John!) arrived at the `quaquaversal' substitution of a
triangular prism in Fig.\ ~\ref{Quaquaversal_substitution}. (I saw no fundamental reason one could
not extend the pinwheel technique to produce matching rules for this
three dimensional example so I never tried, and of course John wasn't
interested in this anyway.)  It was now harder to `see' what tilings
of space are produced, but Fig.\ ~\ref{Quaquaversal_tiling} gives some sense of it.

This was our first serious joint project. Our objective was to
prove various properties of quaquaversal tilings that were not
obtainable in planar tilings such as the pinwheel. For instance if one
looks in a planar substitution tiling the number of orientations of
tiles that appear in an expanding window will grow at most
logarithmically with the diameter of the window, as a simple
consequence of commutativity of rotations.  I hoped one could get
power growth in three dimensions, and in fact we proved this in \cite{CR}. It was
not so easy, and required analyzing the subgroup of the rotation group
through which the tile is turned in the substitution process,
namely the groups generated by rotation about orthogonal axes, one by
$2\pi/3$ and the other by $2\pi/4$. (John liked to represent rotations
by quaternions which was a new world for me.)  We called this $G(3,4)$
and determined a presentation of it and its subgroup $G(3,3)$. (Only
$G(4,4)$ and $G(2,n)$ are finite groups.)  This led to two more papers
(with a third coauthor, Lorenzo Sadun) \cite{CRS1, CRS2} studying groups generated when
the angle between the rotation axes was more interesting geometrically: `geodetic'
angles, whose squared trigonometric functions are rational. (John had
disapproved of my choice of the `pinwheel' name -- it seems the toy I
based it on is called something else in England -- and once we were
working together it was unquestioned that John would produce terms
such as `quaquaversal' and `geodetic'.)

Our results have had a number of applications, and in particular the
quaquaversal tiling has lived up to its original motivation. Both the
pinwheel tilings and the quaquaversal tilings provide easily
computable partitions of (2 and 3 dimensional) space which are
statistically `round': the number of orientations of its components
becomes uniform as the diameter of a window of it grows. They have
therefore been used together with numerical solutions of differential
equations, and other analyses of discrete data in space, for which one
wants to avoid artificial consequences of the symmetries of cubic
partitions. In this regard the pinwheel suffers from the slow
(logarithmic) speed at which it exhibits its roundness and the
quaquaversal successfully overcame this weakness. This subject is
directly related to `expander graphs' and the quaquaversal tilings
have recently been analyzed in that context.

Quaquaversal tilings are produced by successive rotations about
orthogonal axes. To connect better with significant polyhedra, in our
following papers with Sadun we replaced right angles by geodetic
angles, and analyzed the subgroups of $SO(3)$ this produced.  To
understand how we applied this it is useful to step back a bit.

One way the complexity of the rotation group of three dimensions is
exhibited is in fundamental aspects of volume. In two dimensions a
pair of polygons have equal area if and only if they can each be
decomposed into a common set of polygons, i.e.\ they are
rearrangements of a common set of simple components. This is the
Wallace-Bolyai-Gerwien theorem. But this fails badly in three
dimensions (Hilbert's third problem); two polyhedra are
equidecomposable into polyhedra if and only if they have the same
volume {\em and} have the same `Dehn invariant', a result known as
Sydler's theorem. Because of our choice of the class of angles between
rotation axes we were able to apply our results on tiling groups to
compute Dehn invariants for Archimedian polyhedra.

But one can go further into the fundmentals of volume. The
Banach-Tarski paradox shows that two balls of {\em different} volume
can be equidecomposed into a finite number of more complicated
(nonmeasurable) components! The intrusion of measurability is
significant for our story as we will see from another interesting
aperiodic tiling.

\section{The binary tiling of the hyperbolic plane}

\begin{figure}[H]
\center{\includegraphics[width=7in]{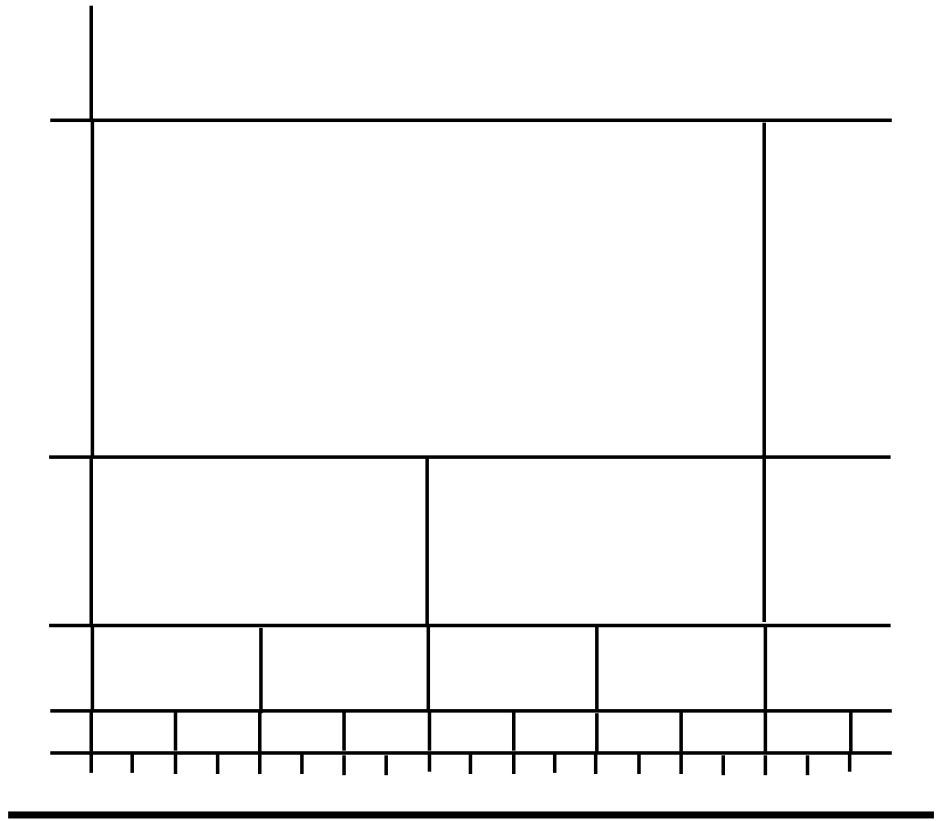}}
\caption{Binary tiling.}
\label{hyperbolic_tiling}
\end{figure}

\begin{figure}[H]
\center{\includegraphics[width=5in]{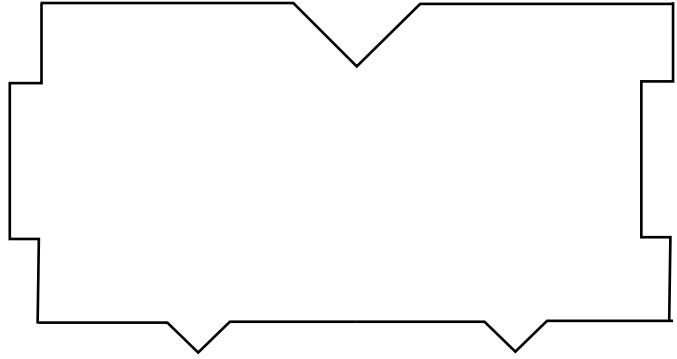}}
\caption{Binary tile.}
\label{hyperbolic_tile}
\end{figure}

In the same paper in which Penrose wrote up the kite and dart tiling
in 1978 \cite{Pen} he also introduced a tiling of the hyperbolic plane,
using congruent copies of a tile two sides of which are geodesics and
two  sides are horocycles. Using the upper half plane model of the
hyperbolic plane, Fig.\ ~\ref{hyperbolic_tiling} shows the tiling 
and Fig.\ ~\ref{hyperbolic_tile} shows the `binary' tile, which has bumps and dents
added. 

In any tiling of the plane made from congruent images of a binary tile, its rectangular
bumps force tiles to uniquely fill out the strip between two
horizontal horocycles, and then the triangular bumps force such strips
to uniquely fill out the half plane below it. There must also be a
strip above it, but it is not unique; there are two ways to place
it. Therefore there are uncountably many ways a given binary tile can sit in
such a tiling, of which Fig.\ ~\ref{hyperbolic_tiling} is representative.

Whether or not the binary tile is `aperiodic' requires a further
reinterpretation of the concept, in terms of symmetries of the hyperbolic plane. It
is complicated but fortunately unnecessary here. We will just note the
connection of these tilings with a famous packing of disks by Karoly
B\"or\"oczky produced a few years earlier \cite{Bor}, in 1974.

\begin{figure}[H]
\center{\includegraphics[width=7in]{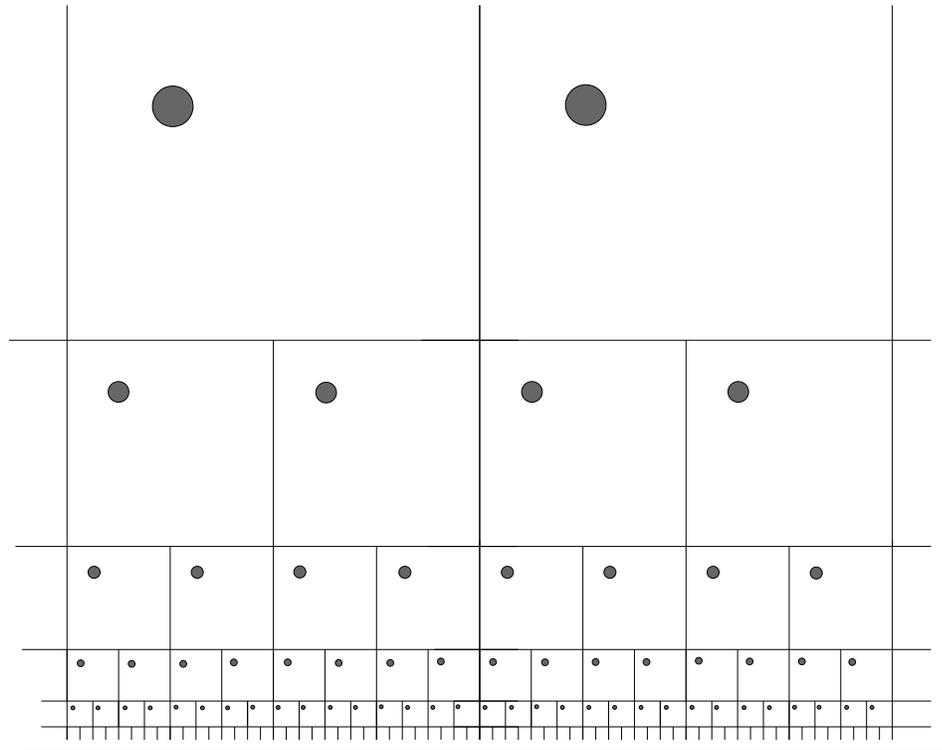}}
\caption{B\"or\"oczky packing.}
\label{Boroczky}
\end{figure}

\begin{figure}[H]
\center{\includegraphics[width=7in]{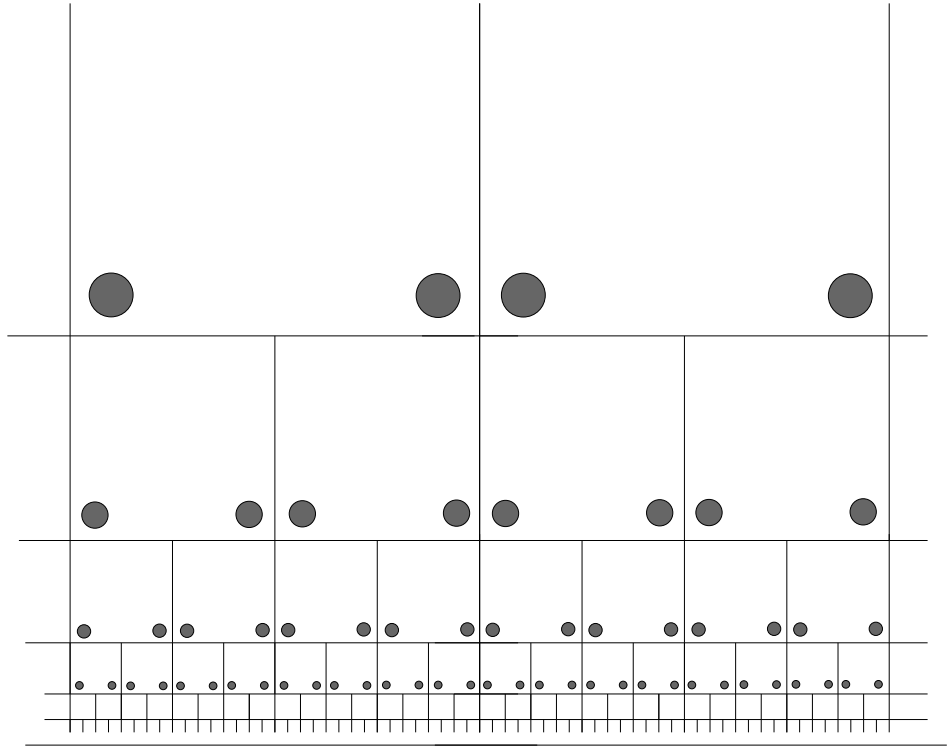}}
\caption{Shifted B\"or\"oczky packing.}
\label{Boroczky2}
\end{figure}

Superimposing a disk in the binary tile one gets disk packings such as
Fig.\ ~\ref{Boroczky}.  Assume the `primary' tile has vertices with complex coordinates
at $ i, 2 + i, 2i$ and $2 + 2i$ and the disk is centered at $1/2 +
7i/4$. Performing the rigid motion $z \to 6z/5$ just on the disks, the primary tile now
has disks at $3/5 + 21i/10$ and $3 + 21i/10$, and the full shifted disk
packing is as in Fig.\ ~\ref{Boroczky2}.

A useful notion of density for the disk packing would require that the
density be unchanged by the rigid motion, but this would contradict
the change from one to two disks in each tile. It took years to come
to grips with this amazing example. 30 years later, an analysis of the
situation by Lewis Bowen in his 2002 thesis\cite{BoR} starting from
aperiodicity and using an ergodic theory perspective, showed that
one can make sense of packing density and aperiodicity in hyperbolic
spaces but one has to avoid pathological or nonmeasureable sets of
packings or tilings.

This is getting beyond the scope of this article, so I'll just sharpen
this last point. The aspect of volume elucidated by
equidecomposability is basically discret mathematics in 2
(Euclidean) dimensions. The Banach-Tarski paradox shows that such a
fundamental study of volume becomes much more complicated in 3
dimensions, involving measurability. Considering volume at a grander
scale, namely the relative volumes of asymptotically large sets, one
arrives at one of the key goals of discrete geometry: the study of
the densest packing of simple bodies (for instance spheres). And in
this realm, B\"or\"oczky's jarring example in the hyperbolic plane 
eventually forced a form of measurability. In summary, analysis
forces its way into this fundamental area of discrete geometry,
because of the complicated nature of the isometry group of the
underlying space.

Penrose introduced two aperiodic tilings in 1978, the kite and dart
tilings in the Euclidean plane and the binary tilings of the
hyperbolic plane. These tilings introduced geometry into aperiodicity,
vastly broadening its appeal, and in particular captured John's
interest. New geometric examples then led to reformulations of
aperiodicity, which also brought in ergodic theory, which John
avoided, but nothing could replace the amazing inventiveness that John
brought to the parts of the story that attracted him. This tribute
to John is built around provocative geometric examples, a playground
in which he was unequaled.

\vfill \eject




\vfill \eject

\vfill \eject
\hbox{}


\end{document}